# Bayesian transformation hazard models

## Gousheng Yin[1] and Joseph G. Ibrahim[2]

*M. D. Anderson Cancer Center and University of North Carolina*

**Abstract:** We propose a class of transformation hazard models for right-censored failure time data. It includes the proportional hazards model (Cox) and the additive hazards model (Lin and Ying) as special cases. Due to the requirement of a nonnegative hazard function, multidimensional parameter constraints must be imposed in the model formulation. In the Bayesian paradigm, the nonlinear parameter constraint introduces many new computational challenges. We propose a prior through a conditional-marginal specification, in which the conditional distribution is univariate, and absorbs all of the nonlinear parameter constraints. The marginal part of the prior specification is free of any constraints. This class of prior distributions allows us to easily compute the full conditionals needed for Gibbs sampling, and hence implement the Markov chain Monte Carlo algorithm in a relatively straightforward fashion. Model comparison is based on the conditional predictive ordinate and the deviance information criterion. This new class of models is illustrated with a simulation study and a real dataset from a melanoma clinical trial.

## 1. Introduction

In survival analysis and clinical trials, the Cox [10] proportional hazards model has been routinely used. For a subject with a possibly time-dependent covariate vector $\mathbf{Z}(t)$, the proportional hazards model is given by,

$$\lambda(t|\mathbf{Z}) = \lambda_0(t) \exp\{\boldsymbol{\beta}'\mathbf{Z}(t)\}, \tag{1.1}$$

where $\lambda_0(t)$ is the unknown baseline hazard function and $\boldsymbol{\beta}$ is the $p \times 1$ parameter vector of interest. Cox [11] proposed to estimate $\boldsymbol{\beta}$ under model (1.1) by maximizing the partial likelihood function and its large sample theory was established by Andersen and Gill [1]. However, the proportionality of hazards might not be a valid modeling assumption in many situations. For example, the true relationship between hazards could be parallel, which leads to the additive hazards model (Lin and Ying [24]),

$$\lambda(t|\mathbf{Z}) = \lambda_0(t) + \boldsymbol{\beta}'\mathbf{Z}(t). \tag{1.2}$$

As opposed to the hazard ratio yielded in (1.1), the hazard difference can be obtained from (1.2), which formulates a direct association between the expected num-ber of events or death occurrences and risk exposures. O'Neill [28] showed that use of the Cox model can result in serious bias when the additive hazards model is correct. Both the multiplicative and additive hazards models have sound biological motivations and solid statistical bases.

[1]Department of Biostatistics & Applied Mathematics, M. D. Anderson Cancer Center, The University of Texas, 1515 Holcombe Boulevard 447, Houston, TX 77030, USA, e-mail: gsyin@mdanderson.org
[2]Department of Biostatistics, The University of North Carolina, Chapel Hill, NC 27599, USA, e-mail: ibrahim@bios.unc.edu

Lin and Ying [25], Martinussen and Scheike [26] and Scheike and Zhang [30] proposed general additive-multiplicative hazards models in which some covariates impose the proportional hazards structure and others induce an additive effect on the hazards. In contrast, we link the additive and multiplicative hazards models in a completely different fashion. Through a simple transformation, we construct a class of hazard-based regression models that includes those two commonly used modeling schemes. In the usual linear regression model, the Box–Cox transformation [4] may be applied to the response variable,

$$(1.3) \qquad \phi(Y) = \begin{cases} (Y^\gamma - 1)/\gamma & \gamma \neq 0 \\ \log(Y) & \gamma = 0, \end{cases}$$

where $\lim_{\gamma \to 0}(Y^\gamma - 1)/\gamma = \log(Y)$. This transformation has been used in survival analysis as well [2, 3, 5, 7, 13, 32]. Breslow and Storer [7] and Barlow [3] applied this family of power transformations to the covariate structure to model the relative risk $R(\mathbf{Z})$,

$$\log R(\mathbf{Z}) = \begin{cases} \{(1 + \boldsymbol{\beta}'\mathbf{Z})^\gamma - 1\}/\gamma & \gamma \neq 0 \\ \log(1 + \boldsymbol{\beta}'\mathbf{Z}) & \gamma = 0, \end{cases}$$

where $R(\mathbf{Z})$ is the ratio of the incidence rate at one level of the risk factor to that at another level. Aranda-Ordaz [2] and Breslow [5] proposed a compromise between these two special cases, $\gamma = 0$ or 1, while their focus was only on grouped survival data by analyzing sequences of contingency tables. Sakia [29] gave an excellent review on this power transformation.

The proportional and additive hazards models may be viewed as two extremes of a family of regression models. On a basis that is very different from the available methods in the literature, we propose a class of regression models for survival data by imposing the Box–Cox transformation on both the baseline hazard $\lambda_0(t)$ and the hazard $\lambda(t|\mathbf{Z})$. This family of transformation models is very general, which includes the Cox proportional hazards model and the additive hazards model as special cases. By adding a transformation parameter, the proposed modeling structure allows a broad class of hazard patterns. In many applications where the hazards are neither proportional nor parallel, our proposed transformation model provides a unified and flexible methodology for analyzing survival data.

The rest of this article is organized as follows. In Section 2.1, we introduce notation and a class of regression models based on the Box–Cox transformed hazards. In Section 2.2, we derive the likelihood function for the proposed model using piecewise constant hazards. In Section 2.3, we propose a prior specification scheme incorporating the parameter constraints within the Bayesian paradigm. In Section 3, we derive the full conditional distributions needed for Gibbs sampling. In Section 4, we introduce model selection methods based on the conditional predictive ordinate (CPO) in Geisser [14] and the deviance information criterion (DIC) proposed by Spiegelhalter *et al.* [31]. We illustrate the proposed methods with data from a melanoma clinical trial, and examine the model using a simulation study in Section 5. We give a brief discussion in Section 6.

## 2. Transformation hazard models

### 2.1. A new class of models

For $n$ independent subjects, let $T_i$ $(i = 1, \ldots, n)$ be the failure time for subject $i$ and $\mathbf{Z}_i(t)$ be the corresponding $p \times 1$ covariate vector. Let $C_i$ be the censoring variable and define $Y_i = \min(T_i, C_i)$. The censoring indicator is $\nu_i = I(T_i \leq C_i)$, where $I(\cdot)$ is the indicator function. Assume that $T_i$ and $C_i$ are independent conditional on $\mathbf{Z}_i(t)$, and that the triplets $\{(T_i, C_i, \mathbf{Z}_i(t)), i = 1, \ldots, n\}$ are independent and identically distributed.

For right-censored failure time data, we propose a class of Box–Cox transformation hazard models,

$$(2.1) \qquad \phi\{\lambda(t|\mathbf{Z}_i)\} = \phi\{\lambda_0(t)\} + \boldsymbol{\beta}'\mathbf{Z}_i(t),$$

where $\phi(\cdot)$ is a known link function given by (1.3). We take $\gamma$ as fixed throughout our development for the following reasons. First, our main goal is to model selection on $\gamma$, by fitting separate models for each value of $\gamma$ and evaluating them through a model selection criterion. Once the best $\gamma$ is chosen according to a model selection criterion, posterior inference regarding $(\boldsymbol{\beta}, \boldsymbol{\lambda})$ is then based on that $\gamma$. Second, in real data settings, there is typically very little information contained in the data to estimate $\gamma$ directly. Third, posterior estimation of $\gamma$ is computationally difficult and often numerically unstable due to the constraint (2.3) as well as its weak identifiability property. To understand how the hazard varies with respect to $\gamma$, we carried out a numerical study as follows. We assume that $\lambda_0(t) = t/3$ in one case, and $\lambda_0(t) = t^2/5$ in another case. A single covariate $Z$ takes a value of 0 or 1 with probability .5, and $\gamma = (0, .25, .5, .75, 1)$. Model (2.1) can be written as

$$\lambda(t|\mathbf{Z}_i) = \{\lambda_0(t)^\gamma + \gamma\boldsymbol{\beta}'\mathbf{Z}_i(t)\}^{1/\gamma}.$$

As shown in Figure 1, there is a broad family of models for $0 \leq \gamma \leq 1$. Our primary interest for $\gamma$ lies in $[0, 1]$, which covers the two popular cases and a family of intermediate modeling structures between the proportional ($\gamma = 0$) and the additive ($\gamma = 1$) hazards models.

Misspecified models may lead to severe bias and wrong statistical inference. In many applications where neither the proportional nor the parallel hazards assumption holds, one can apply (2.1) to the data with a set of prespecified $\gamma$'s, and choose the best fitting model according to a suitable model selection criterion. The need for the general class of models in (2.1) can be demonstrated by the E1690 data from the Eastern Cooperative Oncology Group (ECOG) phase III melanoma clinical trial (Kirkwood *et al.* [23]). The objective of this trial was to compare high-dose interferon to observation (control). Relapse-free survival was a primary outcome variable, which was defined as the time from randomization to progression of tumor or death. As shown in Section 5, the best choice of $\gamma$ in the E1690 data is indeed neither 0 nor 1, but $\gamma = .5$.

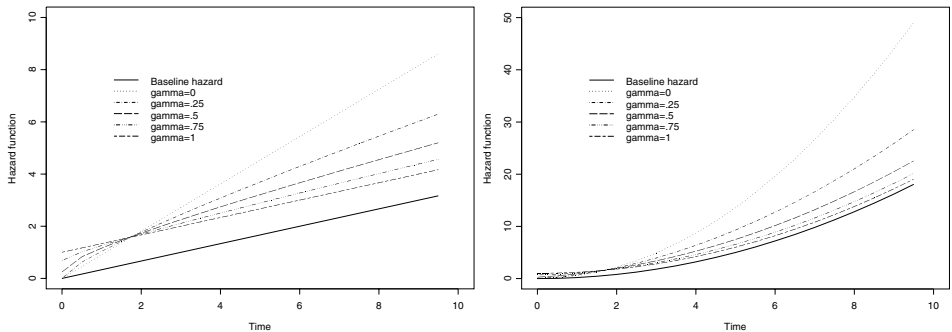

FIG 1. *The relationships between* $\lambda_0(t)$ *and* $\lambda(t|Z) = \{\lambda_0(t)^\gamma + \gamma Z\}^{1/\gamma}$, *with* $Z = 0, 1$. *Left:* $\lambda_0(t) = t/3$; *right:* $\lambda_0(t) = t^2/5$.

Due to the extra parameter $\gamma$, $\boldsymbol{\beta}$ is intertwined with $\lambda_0(t)$ in (2.1). As a result, the model is very different from either the proportional hazards model, which can be solved through the partial likelihood procedure, or the additive hazards model, where the estimating equation can be constructed based on martingale integrals. Here, we propose to conduct inference with this transformation model using a Bayesian approach.

### 2.2. Likelihood function

The piecewise exponential model is chosen for $\lambda_0(t)$. This is a flexible and commonly used modeling scheme and usually serves as a benchmark for the comparison of parametric and nonparametric approaches (Ibrahim, Chen and Sinha [21]). Other nonparametric Bayesian methods for modeling $\lambda_0(t)$ are available in the literature [20, 22, 27]. Let $y_i$ be the observed time for the $i$th subject, $\mathbf{y} = (y_1, \ldots, y_n)'$, $\boldsymbol{\nu} = (\nu_1, \ldots, \nu_n)'$, and $\mathbf{Z}(t) = (\mathbf{Z}_1(t), \ldots, \mathbf{Z}_n(t))'$. Let $J$ denote the number of partitions of the time axis, i.e. $0 < s_1 < \cdots < s_J$, $s_J > y_i$ for $i = 1, \ldots, n$, and that $\lambda_0(y) = \lambda_j$ for $y \in (s_{j-1}, s_j]$, $j = 1, \ldots, J$. When $J = 1$, the model reduces to a parametric exponential model. By increasing $J$, the piecewise constant hazard formulation can essentially model any shape of the underlying hazard. The usual way to partition the time axis is to obtain an approximately equal number of failures in each interval, and to guarantee that each time interval contains at least one failure. Define $\delta_{ij} = 1$ if the $i$th subject fails or is censored in the $j$th interval, and 0 otherwise. Let $D = (n, \mathbf{y}, \mathbf{Z}(t), \boldsymbol{\nu})$ denote the observed data, and $\boldsymbol{\lambda} = (\lambda_1, \ldots, \lambda_J)'$. For ease of exposition and computation, let $\mathbf{Z}_i \equiv \mathbf{Z}_i(t)$, then the likelihood function is

$$
\begin{aligned}
L(\boldsymbol{\beta}, \boldsymbol{\lambda}|D) \;=\; & \prod_{i=1}^{n} \prod_{j=1}^{J} (\lambda_j^\gamma + \gamma \boldsymbol{\beta}' \mathbf{Z}_i)^{\delta_{ij} \nu_i / \gamma} \\
& \times e^{-\delta_{ij}\{(\lambda_j^\gamma + \gamma \boldsymbol{\beta}' \mathbf{Z}_i)^{1/\gamma}(y_i - s_{j-1}) + \sum_{g=1}^{j-1}(\lambda_g^\gamma + \gamma \boldsymbol{\beta}' \mathbf{Z}_i)^{1/\gamma}(s_g - s_{g-1})\}}.
\end{aligned}
\tag{2.2}
$$

### 2.3. Prior distributions

The joint prior distribution of $(\boldsymbol{\beta}, \boldsymbol{\lambda})$ needs to accommodate the nonnegativity constraint for the hazard function, that is,

$$\lambda_j^\gamma + \gamma \boldsymbol{\beta}' \mathbf{Z}_i \geq 0 \quad (i = 1, \ldots, n; j = 1, \ldots, J). \tag{2.3}$$

Constrained parameter problems typically make Bayesian computation and analysis quite complicated [8, 9, 16]. For example, the order constraint on a set of parameters (e.g., $\theta_1 \leq \theta_2 \leq \cdots$) is very common in Bayesian hierarchical models. In these settings, closed form expressions for the normalizing constants in the full conditional distributions are typically available. However, for our model, this is not the case; the normalizing constant involves a complicated intractable integral. The nonnegativity of the hazard constraint is very different from the usual order constraints. If the hazard is negative, the likelihood function and the posterior density are not well defined. One way to proceed with this nonlinear constraint is to specify an appropriately truncated joint prior distribution for $(\boldsymbol{\beta}, \boldsymbol{\lambda})$, such as a truncated multivariate normal prior $N(\boldsymbol{\mu}, \boldsymbol{\Sigma})$ for $(\boldsymbol{\beta}|\boldsymbol{\lambda})$ to satisfy this constraint. This would lead to a prior distribution of the form

$$\pi(\boldsymbol{\beta}, \boldsymbol{\lambda}) = \pi(\boldsymbol{\beta}|\boldsymbol{\lambda})\pi(\boldsymbol{\lambda})I(\lambda_j^\gamma + \gamma\boldsymbol{\beta}'\mathbf{Z}_i \geq 0,\ i=1,\ldots,n;\ j=1,\ldots,J).$$

Following this route, we would need to analytically compute the normalizing constant,

$$c(\boldsymbol{\lambda}) = \int \cdots \int_{\lambda_j^\gamma + \gamma\boldsymbol{\beta}'\mathbf{Z}_i \geq 0 \text{ for all } i,j} \exp\left\{-\frac{1}{2}(\boldsymbol{\beta}-\boldsymbol{\mu})'\boldsymbol{\Sigma}^{-1}(\boldsymbol{\beta}-\boldsymbol{\mu})\right\} d\beta_1 \cdots d\beta_p$$

to construct the full conditional distribution of $\boldsymbol{\lambda}$. However, $c(\boldsymbol{\lambda})$ involves a $p$-dimensional integral on a complex nonlinear constrained parameter space, which cannot be obtained in a closed form. Such a prior would lead to intractable full conditionals, therefore making Gibbs sampling essentially impossible.

To circumvent the multivariate constrained parameter problem, we reduce our prior specification to a one-dimensional truncated distribution, and thus the normalizing constant can be obtained in a closed form. Without loss of generality, we assume that all the covariates are positive. Let $\mathbf{Z}_{i(-k)}$ denote the covariate $\mathbf{Z}_i$ with the $k$th component $Z_{ik}$ deleted, and let $\boldsymbol{\beta}_{(-k)}$ denote the $(p-1)$-dimensional parameter vector with $\beta_k$ removed, and define

$$h_\gamma(\lambda_j, \boldsymbol{\beta}_{(-k)}, \mathbf{Z}_i) = \min_{i,j}\left\{\frac{\lambda_j^\gamma + \gamma\boldsymbol{\beta}'_{(-k)}\mathbf{Z}_{i(-k)}}{\gamma Z_{ik}}\right\}.$$

We propose a joint prior for $(\boldsymbol{\beta}, \boldsymbol{\lambda})$ of the form

$$(2.4) \qquad \pi(\boldsymbol{\beta}, \boldsymbol{\lambda}) = \pi(\beta_k|\boldsymbol{\beta}_{(-k)}, \boldsymbol{\lambda})I\left(\beta_k \geq -h_\gamma(\lambda_j, \boldsymbol{\beta}_{(-k)}, \mathbf{Z}_i)\right)\pi(\boldsymbol{\beta}_{(-k)}, \boldsymbol{\lambda}).$$

We see that $\beta_k$ and $(\boldsymbol{\beta}_{(-k)}, \boldsymbol{\lambda})$ are not independent a priori due to the nonlinear parameter constraint. This joint prior specification only involves one parameter $\beta_k$ in the constraints and makes all the other parameters $(\boldsymbol{\beta}_{(-k)}, \boldsymbol{\lambda})$ free of constraints.

Let $\Phi(\cdot)$ denote the cumulative distribution function of the standard normal distribution. Specifically, we take $(\beta_k|\boldsymbol{\beta}_{(-k)}, \boldsymbol{\lambda})$ to have a truncated normal distribution,

$$(2.5) \qquad \pi(\beta_k|\boldsymbol{\beta}_{(-k)}, \boldsymbol{\lambda}) = \frac{\exp\{-\frac{\beta_k^2}{2\sigma_k^2}\}}{c(\boldsymbol{\beta}_{(-k)}, \boldsymbol{\lambda})}I\left(\beta_k \geq -h_\gamma(\lambda_j, \boldsymbol{\beta}_{(-k)}, \mathbf{Z}_i)\right),$$

where the normalizing constant depends on $\boldsymbol{\beta}_{(-k)}$ and $\boldsymbol{\lambda}$, given by

$$(2.6) \qquad c(\boldsymbol{\beta}_{(-k)}, \boldsymbol{\lambda}) = \sqrt{2\pi}\sigma_k\left[1 - \Phi\left(-\frac{h_\gamma(\lambda_j, \boldsymbol{\beta}_{(-k)}, \mathbf{Z}_i)}{\sigma_k}\right)\right].$$

Thus, we need only to constrain one parameter $\beta_k$ to guarantee the nonnegativity of the hazard function and allow the other parameters, $(\boldsymbol{\beta}_{(-k)}, \boldsymbol{\lambda})$, to be free.

Although not required for the development, we can take $\boldsymbol{\beta}_{(-k)}$ and $\boldsymbol{\lambda}$ to be independent a priori in (2.4), $\pi(\boldsymbol{\beta}_{(-k)}, \boldsymbol{\lambda}) = \pi(\boldsymbol{\beta}_{(-k)})\pi(\boldsymbol{\lambda})$. In addition, we can specify a normal prior distribution for each component of $\boldsymbol{\beta}_{(-k)}$. We assume that the components of $\boldsymbol{\lambda}$ are independent a priori, and each $\lambda_j$ has a Gamma$(\alpha, \xi)$ distribution.

## 3. Gibbs sampling

For $0 \leq \gamma \leq 1$, it can be shown that the full conditionals of $(\beta_1, \ldots, \beta_p)$ are log-concave, in which case we only need to use the adaptive rejection sampling (ARS) algorithm proposed by Gilks and Wild [19]. Due to the non-log-concavity of the full conditionals of the $\lambda_j$'s, a Metropolis step is required within the Gibbs steps, for details see Gilks, Best and Tan [18]. For each Gibbs sampling step, the support for the parameter to be sampled is set to satisfy the constraint (2.3), such that the likelihood function is well defined within the sampling range. For $i = 1, \ldots, n; j = 1, \ldots, J; k = 1, \ldots, p$, the following inequalities need to be satisfied,

$$\beta_k \geq -h_\gamma(\lambda_j, \boldsymbol{\beta}_{(-k)}, \mathbf{Z}_i), \quad \lambda_j \geq -\min_i\{(\gamma\boldsymbol{\beta}'\mathbf{Z}_i)^{1/\gamma}, 0\}.$$

Suppose that the $k$th component of $\boldsymbol{\beta}$ has a truncated normal prior as given in (2.5), and all other parameters are left free. The full conditionals of the parameters are given as follows:

$$\begin{aligned}
\pi(\beta_k|\boldsymbol{\beta}_{(-k)}, \boldsymbol{\lambda}, D) &\propto L(\boldsymbol{\beta}, \boldsymbol{\lambda}|D)\pi(\beta_k|\boldsymbol{\beta}_{(-k)}, \boldsymbol{\lambda}) \\
\pi(\beta_l|\boldsymbol{\beta}_{(-l)}, \boldsymbol{\lambda}, D) &\propto L(\boldsymbol{\beta}, \boldsymbol{\lambda}|D)\pi(\beta_l)/c(\boldsymbol{\beta}_{(-k)}, \boldsymbol{\lambda}) \\
\pi(\lambda_j|\boldsymbol{\beta}, \boldsymbol{\lambda}_{(-j)}, D) &\propto L(\boldsymbol{\beta}, \boldsymbol{\lambda}|D)\pi(\lambda_j)/c(\boldsymbol{\beta}_{(-k)}, \boldsymbol{\lambda})
\end{aligned}$$

where

$$\begin{aligned}
\pi(\beta_l) &\propto \exp\{-\beta_l^2/(2\sigma_l^2)\}, \quad l \neq k, l = 1, \ldots, p, \\
\pi(\lambda_j) &\propto \lambda_j^{\alpha-1}\exp(-\xi\lambda_j), \quad j = 1, \ldots, J.
\end{aligned}$$

These full conditionals have nice tractable structures, since $c(\boldsymbol{\beta}_{(-k)}, \boldsymbol{\lambda})$ has a closed form with our proposed prior specification. Posterior estimation is very robust with respect to the conditioning scheme (the choice of $k$) in (2.4).

## 4. Model assessment

It is crucial to compare a class of competing models for a given dataset and select the model that best fits the data. After fitting the proposed models for a set of pre-specified $\gamma$'s, we compute the CPO and DIC statistics, which are the two commonly used measures of model adequacy [14, 15, 12, 31].

We first introduce the CPO as follows. Let $\mathbf{Z}^{(-i)}$ denote the $(n-1) \times p$ covariate matrix with the $i$th row deleted, let $\mathbf{y}^{(-i)}$ denote the $(n-1) \times 1$ response vector with $y_i$ deleted, and $\boldsymbol{\nu}^{(-i)}$ is defined similarly. The resulting data with the $i$th case deleted can be written as $D^{(-i)} = \{(n-1), \mathbf{y}^{(-i)}, \mathbf{Z}^{(-i)}, \boldsymbol{\nu}^{(-i)}\}$. Let $f(y_i|\mathbf{Z}_i, \boldsymbol{\beta}, \boldsymbol{\lambda})$ denote the density function of $y_i$, and let $\pi(\boldsymbol{\beta}, \boldsymbol{\lambda}|D^{(-i)})$ denote the posterior density of $(\boldsymbol{\beta}, \boldsymbol{\lambda})$ given $D^{(-i)}$. Then, CPO$_i$ is the marginal posterior predictive density of

$y_i$ given $D^{(-i)}$, which can be written as

$$\begin{aligned}
\text{CPO}_i &= f(y_i|\mathbf{Z}_i, D^{(-i)}) \\
&= \int\int f(y_i|\mathbf{Z}_i, \boldsymbol{\beta}, \boldsymbol{\lambda})\pi(\boldsymbol{\beta}, \boldsymbol{\lambda}|D^{(-i)})d\boldsymbol{\beta}d\boldsymbol{\lambda} \\
&= \left\{\int\int \frac{\pi(\boldsymbol{\beta}, \boldsymbol{\lambda}|D)}{f(y_i|\mathbf{Z}_i, \boldsymbol{\beta}, \boldsymbol{\lambda})}d\boldsymbol{\beta}d\boldsymbol{\lambda}\right\}^{-1}.
\end{aligned}$$

For the proposed transformation model, a Monte Carlo approximation of $\text{CPO}_i$ is given by,

$$\widehat{\text{CPO}}_i = \left\{\frac{1}{M}\sum_{m=1}^{M} \frac{1}{L_i(\boldsymbol{\beta}_{[m]}, \boldsymbol{\lambda}_{[m]}|y_i, \mathbf{Z}_i, \nu_i)}\right\}^{-1},$$

where

$$\begin{aligned}
L_i(\boldsymbol{\beta}_{[m]}, \boldsymbol{\lambda}_{[m]}|y_i, \mathbf{Z}_i, \nu_i) &= \prod_{j=1}^{J}(\lambda_{j,[m]}^{\gamma} + \gamma\boldsymbol{\beta}'_{[m]}\mathbf{Z}_i)^{\delta_{ij}\nu_i/\gamma} \\
&\quad \times \exp\left[-\delta_{ij}\left\{(\lambda_{j,[m]}^{\gamma} + \gamma\boldsymbol{\beta}'_{[m]}\mathbf{Z}_i)^{1/\gamma}(y_i - s_{j-1})\right.\right. \\
&\quad\left.\left. + \sum_{g=1}^{j-1}(\lambda_{g,[m]}^{\gamma} + \gamma\boldsymbol{\beta}'_{[m]}\mathbf{Z}_i)^{1/\gamma}(s_g - s_{g-1})\right\}\right].
\end{aligned}$$

Note that $M$ is the number of Gibbs samples after burn-in, and $\boldsymbol{\lambda}_{[m]} = (\lambda_{1,[m]}, \ldots, \lambda_{J,[m]})'$ and $\boldsymbol{\beta}_{[m]}$ are the samples of the $m$th Gibbs iteration. A common summary statistic based on the $\text{CPO}_i$'s is $B = \sum_{i=1}^{n}\log(\text{CPO}_i)$, which is often called the logarithm of the pseudo Bayes factor. A larger value of $B$ indicates a better fit of a model.

Another model assessment criterion is the DIC (Spiegelhalter et al. [31]), defined as

$$\text{DIC} = 2\overline{\text{Dev}(\boldsymbol{\beta}, \boldsymbol{\lambda})} - \text{Dev}(\bar{\boldsymbol{\beta}}, \bar{\boldsymbol{\lambda}}),$$

where $\text{Dev}(\boldsymbol{\beta}, \boldsymbol{\lambda}) = -2\log L(\boldsymbol{\beta}, \boldsymbol{\lambda}|D)$ is the deviance, and $\overline{\text{Dev}(\boldsymbol{\beta}, \boldsymbol{\lambda})}$, $\bar{\boldsymbol{\beta}}$ and $\bar{\boldsymbol{\lambda}}$ are the corresponding posterior means. Specifically, in our proposed model,

$$\text{DIC} = -\frac{4}{M}\sum_{m=1}^{M}\log L(\boldsymbol{\beta}_{[m]}, \boldsymbol{\lambda}_{[m]}|D) + 2\log L(\bar{\boldsymbol{\beta}}, \bar{\boldsymbol{\lambda}}|D).$$

The smaller the DIC value, the better the fit of the model.

## 5. Numerical studies

### 5.1. Application

As an illustration, we applied the transformation models to the E1690 data. There were a total of $n = 427$ patients on these combined treatment arms. The covariates in this analysis were treatment (high-dose interferon or observation), age (a continuous variable which ranged from 19.13 to 78.05 with mean 47.93 years), sex (male or female) and nodal category (1 if there were no positive nodes, or 2 otherwise). Figure 2 shows the estimated cumulative hazard curves for the interferon and observation groups based on the Nelson–Aalen estimator.

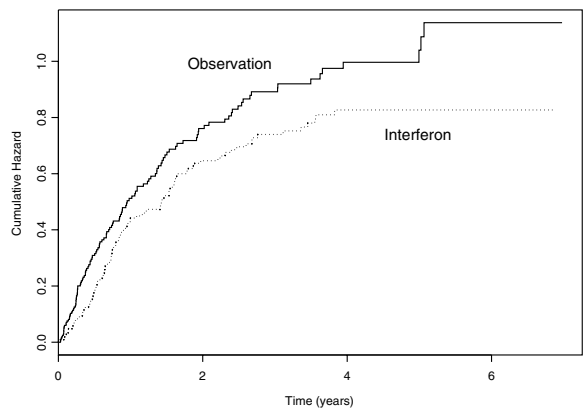

Fig 2. *The estimated cumulative hazard curves for the two arms in E1690*

Table 1
*The B/DIC statistics with respect to $\gamma$ and $J$ in the E1690 data*

|  |  | $J$ | | |
|---|---|---|---|---|
|  |  | 1 | 5 | 10 |
| | 0 | $-\mathbf{567.43/1129.19}$ | $-528.36/1051.84$ | $-555.46/1105.48$ |
| | .25 | $-567.96/1131.71$ | $-523.74/1045.68$ | $-534.57/1066.86$ |
| $\gamma$ | .5 | $-568.47/1133.72$ | $-\mathbf{522.55/1043.64}$ | $-529.13/1056.44$ |
| | .75 | $-568.89/1135.16$ | $-522.66/1043.86$ | $-527.47/1053.17$ |
| | 1 | $-569.46/1136.54$ | $-523.04/1044.84$ | $-\mathbf{526.80/1052.06}$ |

We constrained the regression coefficient for treatment, $\beta_1$, to have the truncated normal prior. We prespecified $\gamma = (0, .25, .5, .75, 1)$ and took the priors for $\boldsymbol{\beta} = (\beta_1, \beta_2, \beta_3, \beta_4)'$ and $\boldsymbol{\lambda} = (\lambda_1, \ldots, \lambda_J)'$ to be noninformative. For example, $(\beta_1 | \boldsymbol{\lambda}, \boldsymbol{\beta}_{(-1)})$ was assigned the truncated $N(0, 10,000)$ prior as defined in (2.5), $(\beta_l, l = 2, 3, 4)$ were taken to have independent $N(0, 10,000)$ prior distributions, and $\lambda_j \sim \text{Gamma}(2, .01)$, and independent for $j = 1, \ldots, J$. To allow for a fair comparison between different models using different $\gamma$'s, we used the same noninformative priors across all the targeted models.

The shape of the baseline hazard function is controlled by $J$. The finer the partition of the time axis, the more general the pattern of the hazard function that is captured. However, by increasing $J$, we introduce more unknown parameters (the $\lambda_j$'s). For the proposed transformation model, $\gamma$ also directly affects the shape of the hazard function, and specifically, there is much interplay between $J$ and $\gamma$ in controlling the shape of the hazard, and in some sense $\gamma$ and $J$ are somewhat confounded. Thus when searching for the best fitting model, we must find suitable $J$ and $\gamma$ simultaneously. Similar to a grid search, we set $J = (1, 5, 10)$, and located the point $(J, \gamma)$ that yielded the largest $B$ statistic and the smallest DIC.

After a burn-in of 2,000 samples and thinned by 5 iterations, the posterior computations were based on 10,000 Gibbs samples. The $B$ and DIC statistics for model selection are summarized in Table 1. The two model selection criteria are quite consistent with each other, and both lead to the same best model with $J = 5$ and $\gamma = .5$. Table 2 summarizes the posterior means, standard deviations and the 95% highest posterior density (HPD) intervals for $\boldsymbol{\beta}$ using $J = (1, 5, 10)$ and $\gamma = (0, .5, 1)$. For the best model (with $J = 5$ and $\gamma = .5$), we see that the treatment effect has a 95% HPD interval that does not include 0, confirming that treatment with high-dose

TABLE 2
*Posterior means, standard deviations, and 95% HPD intervals for the E1690 data*

| $J$ | $\gamma$ | Covariate | Mean | SD | 95% HPD Interval |
|---|---|---|---|---|---|
| 1 | 0 | Treatment | −.2888 | .1299 | (−.5369, −.0310) |
|   |   | Age | .0117 | .0050 | (.0016, .0214) |
|   |   | Sex | −.3479 | .1375 | (−.6372, −.0962) |
|   |   | Nodal Category | .5267 | .1541 | (.2339, .8346) |
|   | .5 | Treatment | −.1398 | .0626 | (−.2588, −.0111) |
|   |   | Age | .0056 | .0024 | (.0011, .0103) |
|   |   | Sex | −.1464 | .0644 | (−.2791, −.0254) |
|   |   | Nodal Category | .2179 | .0688 | (.0835, .3529) |
|   | 1 | Treatment | −.0655 | .0299 | (−.1245, −.0078) |
|   |   | Age | .0026 | .0011 | (.0004, .0047) |
|   |   | Sex | −.0593 | .0293 | (−.1155, −.0007) |
|   |   | Nodal Category | .0863 | .0296 | (.0304, .1471) |
| 5 | 0 | Treatment | −.4865 | .1295 | (−.7492, −.2408) |
|   |   | Age | −.0036 | .0050 | (−.0133, .0061) |
|   |   | Sex | −.4423 | .1421 | (−.7196, −.1684) |
|   |   | Nodal Category | .1461 | .1448 | (−.1307, .4298) |
|   | .5 | Treatment | −.1835 | .0626 | (−.3066, −.0604) |
|   |   | Age | .0017 | .0024 | (−.0030, .0064) |
|   |   | Sex | −.1557 | .0655 | (−.2853, −.0310) |
|   |   | Nodal Category | .1141 | .0685 | (−.0179, .2510) |
|   | 1 | Treatment | −.0525 | .0274 | (−.1058, .0007) |
|   |   | Age | .0011 | .0009 | (−.0006, .0027) |
|   |   | Sex | −.0334 | .0249 | (−.0818, .0148) |
|   |   | Nodal Category | .0265 | .0224 | (−.0169, .0705) |
| 10 | 0 | Treatment | −.7238 | .1260 | (−.9639, −.4710) |
|   |   | Age | −.0175 | .0047 | (−.0269, −.0084) |
|   |   | Sex | −.6368 | .1439 | (−.9158, −.3544) |
|   |   | Nodal Category | .1685 | .1302 | (−.4184, .0859) |
|   | .5 | Treatment | −.2272 | .0629 | (−.3581, −.1094) |
|   |   | Age | −.0009 | .0023 | (−.0056, .0035) |
|   |   | Sex | −.1791 | .0649 | (−.3094, −.0546) |
|   |   | Nodal Category | .0534 | .0670 | (−.0814, .1798) |
|   | 1 | Treatment | −.0610 | .0274 | (−.1142, −.0070) |
|   |   | Age | .0006 | .0008 | (−.0010, .0021) |
|   |   | Sex | −.0334 | .0256 | (−.0850, .0155) |
|   |   | Nodal Category | .0107 | .0225 | (−.0325, .0569) |

interferon indeed substantially reduced the risk of melanoma relapse compared to observation.

In Figure 3, we present the estimated hazards for the interferon and observation arms for $\gamma = 0$, .5 and 1 using $J = 5$. It is important to note that, when $\gamma = .5$, the hazard ratio increases over time while the hazard difference decreases.

The proportional hazards model yields a hazard ratio of 1.63, the additive hazards model gives a hazard difference of .05, and the model with $\gamma = .5$ shows hazard ratios of 1.27, 1.36 and 1.61, and hazard differences of .14, .11 and .07 at .5, 1 and 3 years, respectively. This interesting feature between the hazards cannot be captured through a conventional modeling structure. An opposite phenomenon in which the difference of the hazards increases in $t$ whereas their ratio decreases, was noted in the British doctors study (Breslow and Day [6], p.112, pp. 336-338), which examined the effects of cigarette smoking on mortality. We also computed the half year and one year posterior predictive survival probabilities for a 48 years old male patient under the high-dose interferon treatment with one or more positive nodes. When $\gamma = .5$, the .5 year posterior predictive survival probabilities are .8578, .7686 and .7804 for $J = 1$, 5 and 10; the 1 year survival probabilities are .7357, .6043 and .6240, respectively. When $J$ is large enough, the posterior inference becomes stable.

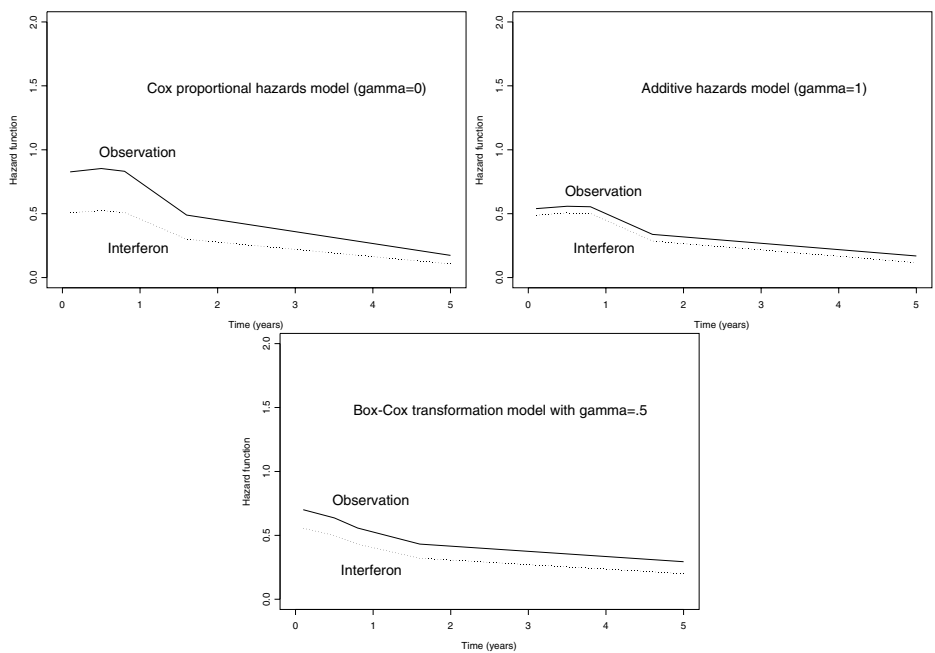

FIG 3. *Estimated hazards under models with $\gamma = 0$, .5 and 1, for male subjects at age = 47.93 years and with one or more positive nodes, using $J = 5$.*

TABLE 3
*Sensitivity analysis with $\beta_k$ having a truncated normal prior using $J = 5$ and $\gamma = .5$*

| Truncated Covariate | Regression Coefficient | Mean | SD | 95% HPD Interval |
|---|---|---|---|---|
| *Age* | Treatment | −.1862 | .0633 | (−.3122, −.0627) |
| | Age | .0016 | .0024 | (−.0032, .0063) |
| | Sex | −.1551 | .0665 | (−.2802, −.0187) |
| | Nodal Category | .1132 | .0697 | (−.0229, .2511) |
| *Sex* | Treatment | −.1883 | .0634 | (−.3107, −.0592) |
| | Age | .0017 | .0024 | (−.0032, .0063) |
| | Sex | −.1572 | .0651 | (−.2801, −.0296) |
| | Nodal Category | .1131 | .0672 | (−.0165, .2448) |
| *Nodal Category* | Treatment | −.1850 | .0633 | (−.3037, −.0566) |
| | Age | .0017 | .0024 | (−.0030, .0062) |
| | Sex | −.1519 | .0662 | (−.2819, −.0236) |
| | Nodal Category | .1124 | .0679 | (−.0223, .2416) |

We examined MCMC convergence based on the method proposed by Geweke [17]. The Markov chains mixed well and converged fast. We conducted a sensitivity analysis on the choice of the conditioning scheme in the prior (2.5) by choosing the regression coefficient of each covariate to have a truncated normal prior. The results in Table 3 show the robustness of the model to the choice of the constrained parameter in the prior specification. This demonstrates the appealing feature of the proposed prior specification, which thus facilitates an attractive computational procedure.

## 5.2. Simulation

We conducted a simulation study to examine properties of the proposed model. The failure times were generated from model (2.1) with $\gamma = .5$. We assumed a constant baseline hazard, i.e., $\lambda_0(t) = .5$, and two covariates were generated independently: $Z_1 \sim N(5, 1)$ and $Z_2$ is a binary random variable taking a value of 1 or 2 with probability .5. The corresponding regression parameters were $\beta_1 = .7$ and $\beta_2 = 1$. The censoring times were simulated from a uniform distribution to achieve approximately a 25% censoring rate. The sample sizes were $n = 300$, 500 and 1,000, and we replicated 500 simulations for each configuration.

TABLE 4
*Simulation results based on 500 replications, with the true values $\beta_1 = .7$ and $\beta_2 = 1$*

| $n$ | $c\%$ | Mean ($\beta_1$) | SD ($\beta_1$) | Mean ($\beta_2$) | SD ($\beta_2$) |
|---|---|---|---|---|---|
| 300 | 0 | .7705 | .2177 | 1.0556 | .4049 |
|  | 25 | .7430 | .2315 | 1.0542 | .4534 |
| 500 | 0 | .7424 | .1989 | 1.0483 | .3486 |
|  | 25 | .7510 | .2084 | 1.0503 | .3781 |
| 1000 | 0 | .7273 | .1784 | 1.0412 | .2920 |
|  | 25 | .7394 | .1869 | 1.0401 | .3100 |

Noninformative prior distributions were specified for the unknown parameters as in the E1690 example. For each Markov chain, we took a burn-in of 200 samples and the posterior estimates were based on 5,000 Gibbs samples. The posterior means and standard deviations are summarized in Table 4, which show the convergence of the posterior means of the parameters to the true values. As the sample size increases, the posterior means of $\beta_1$ and $\beta_2$ approach their true values and the corresponding standard deviations decrease. As the censoring rate increases, the posterior standard deviation also increases.

## 6. Discussion

We have proposed a class of survival models based on the Box–Cox transformed hazard functions. This class of transformation models makes hazard-based regression more flexible, general, and versatile than other methods, and opens a wide family of relationships between the hazards. Due to the complexity of the model, we have proposed a joint prior specification scheme by absorbing the non-linear constraint into one parameter while leaving all the other parameters free of constraints. This prior specification is quite general and can be applied to a much broader class of constrained parameter problems arising from regression models. It is usually difficult to interpret the parameters in the proposed model except when $\gamma = 0$ or 1. However, if the primary aim is for prediction of survival, the best fitting Box–Cox transformation model could be useful.

### Acknowledgements

We would like to thank Professor Javier Rojo and anonymous referees for helpful comments which led to great improvement of the article.